\documentclass[12pt, a4paper]{amsart}
\newtheorem{teo}{Theorem}[section]
\newtheorem{lema}{Lemma}[section]
\newtheorem{coro}{Corollary}[section]
\newtheorem{defi}{Definition}[section]
\newtheorem{rem}{Remark}[section]

\def\ep{\varepsilon}

\def\ve{\varepsilon}
\def\RR{{\mathbb{R}}}

\def\Oc{\RR^N\setminus\Omega}
\def\uep{u^\ep}
\def\wep{w^\ep}
\def\di{\displaystyle}

\begin{document}

\title[Nonlocal diffusion
problems]{How to approximate the heat equation with Neumann
boundary conditions by nonlocal diffusion problems}

\author[C. Cortazar\and M. Elgueta \and J.D. Rossi \and N. Wolanski]{Carmen
Cortazar
\and Manuel Elgueta \and Julio D. Rossi \and Noemi Wolanski}

% 35K57= Reaction-diffusion equations.
% 35B40= Asymptotic behavior of solutions.
\address{Carmen Cortazar and Manuel Elgueta \hfill\break\indent
Departamento  de Matem{\'a}tica, Universidad Catolica de Chile,
\hfill\break\indent Casilla 306, Correo 22, Santiago, Chile. }
\email{\tt ccortaza@mat.puc.cl, melgueta@mat.puc.cl.}

\address{Julio D. Rossi \hfill\break\indent
Consejo Superior de Investigaciones Cient\'{\i}ficas (CSIC),
\hfill\break\indent Serrano 123, Madrid, Spain,
\hfill\break\indent on leave from Departamento de Matem{\'a}tica,
FCEyN \hfill\break\indent UBA (1428) Buenos Aires, Argentina. }
\email{{\tt jrossi@dm.uba.ar}}

\address{Noemi Wolanski \hfill\break\indent
Departamento  de Matem{\'a}tica, FCEyN \hfill\break\indent UBA (1428)
Buenos Aires, Argentina.} \email{{\tt wolanski@dm.uba.ar} }

\date{}

\begin{abstract}
We present a model for nonlocal diffusion with Neumann boundary
conditions in a bounded smooth domain prescribing the flux through
the boundary. We study the limit of this family of nonlocal
diffusion operators when a rescaling parameter related to the
kernel of the nonlocal operator goes to zero. We prove that the
solutions of this family of problems converge to a solution of the
heat equation with Neumann boundary conditions.
\end{abstract}

\maketitle

\date{}

\section{Introduction}
\label{convergence.heat.eq} \setcounter{equation}{0}

The purpose of this article is to show that the solutions of the
usual Neumann boundary value problem for the heat equation can be
approximated by solutions of a sequence of nonlocal ``Neumann''
boundary value problems.

Let $J: \RR^N \to \RR$ be a nonnegative, radial, continuous
function with $\int_{\RR^N} J(z)\, dz =1$. Assume also that $J$ is
strictly positive in $B(0,d)$ and vanishes in $\RR^N \setminus
B(0,d)$. Nonlocal evolution equations of the form
\begin{equation} \label{11}
u_t (x,t) = (J*u-u) (x,t) = \int_{\RR^N} J(x-y)u(y,t) \, dy -
u(x,t) ,
\end{equation}
and variations of it, have been recently widely used to model
diffusion processes. More precisely, as stated in \cite{F}, if
$u(x,t)$ is thought of as a density at the point $x$ at time $t$
and $J(x-y)$ is thought of as the probability distribution of
jumping from location $y$ to location $x$, then $\int_{\RR^N}
J(y-x)u(y,t)\, dy = (J*u)(x,t)$ is the rate at which individuals
are arriving at position $x$ from all other places and $-u(x,t) =
-\int_{\RR^N} J(y-x)u(x,t)\, dy$ is the rate at which they are
leaving location $x$ to travel to all other sites. This
consideration, in the absence of external or internal sources,
leads immediately to  the fact that the density $u$ satisfies
equation (\ref{11}). For recent references on nonlocal diffusion
see, \cite{BCh},
\cite{BCh2},
\cite{BFRW},
\cite{BH}, \cite{BH2},
\cite{CF}, \cite{C},
 \cite{F},
\cite{LW}, \cite{W}, \cite{Z} and references therein.

Given a bounded, connected and smooth domain $\Omega$, one of the
most common boundary conditions that has been imposed in the
literature to the heat equation, $u_t = \Delta u$,  is the {\it
Neumann boundary condition}, $\partial u /
\partial \eta (x,t)=g(x,t)$, $x\in
\partial \Omega$, which leads to the following classical problem,
\begin{equation}\label{calor.I}
\begin{cases}
u_t-\Delta u=0\quad&\mbox{in}\quad \Omega\times(0,T),\\
\di\frac{\partial u}{\partial\eta}=g\quad&\mbox{on}\quad
\partial\Omega\times(0,T),\\
u(x,0)=u_0(x)\quad&\mbox{in}\quad \Omega.
\end{cases}
\end{equation}

In this article we propose a nonlocal ``Neumann'' boundary value
problem, namely
 \begin{equation} \label{Neumann}
 u_t (x,t) =  \displaystyle \int_{\Omega} J(x-y)\big(u(y,t) - u(x,t)\big)
 \, dy   + \displaystyle \int_{\RR^N \setminus \Omega} G(x,x-y)
 g(y,t) \, dy ,
 \end{equation}
where $G(x,\xi)$ is smooth and compactly supported in $\xi$
uniformly in $x$.

In this model the first integral takes into account the diffusion
inside $\Omega$. In fact, as we have explained, the integral $\int
J(x-y) (u(y,t) - u(x,t)) \, dy$ takes into account the individuals
arriving or leaving position $x$ from or to other places. Since we
are integrating in $\Omega$, we are imposing that diffusion takes
place only in $\Omega$. The last term takes into account the
prescribed flux of individuals that enter or leave the domain.

The nonlocal Neumann model
\eqref{Neumann} and the Neumann
problem for the heat equation \eqref{calor.I} share many
properties. For example,  a comparison principle holds for both
equations when $G$ is nonnegative and the asymptotic behavior of
their solutions as $t\to
\infty$ is similar, see \cite{cerw}.

Existence and uniqueness of solutions of
\eqref{Neumann} with general $G$ is proved by a fixed point argument in Section 2.
Also, a comparison principle when $G\ge 0$ is proved in that
section.

Our main goal is to show that the Neumann problem for the heat
equation \eqref{calor.I} can be approximated by suitable nonlocal
Neumann problems
\eqref{Neumann}.

More precisely, for given $J$ and $G$ we consider the rescaled
kernels
\begin{equation}
\label{rescales}
J_\ve (\xi ) =
C_1 \frac{1}{\ve^{N} } J \left(\frac{\xi}{\ve}\right),\qquad
 G_\ve (x, \xi ) =
C_1 \frac{1}{\ve^{N} }\, G \left(x, \frac{\xi}{\ve}\right)
\end{equation}
with
$$
C_1^{-1}=\frac12\int_{B(0,d)}J(z)z_N^2\,dz,
$$
which is a normalizing constant in order to obtain the Laplacian
in the limit instead of a multiple of it. Then, we consider the
solution $\uep (x,t)$ to
\begin{equation} \label{Neumann.C.1}
\left\{
\begin{array}{rl}
\uep_t (x,t) &=  \displaystyle \frac{1}{\ve^2} \int_{\Omega}
 J_\ve( x-y)(\uep(y,t) - \uep(x,t)) \, dy\\
 &\hskip2cm+\displaystyle \frac{1 }{\ve} \int_{\Oc} G_\ve (x, x-y)
g(y,t) \, dy,\\
\uep(x,0)&=u_0(x).
\end{array}
\right.
\end{equation}

We prove in this paper that
$$
u^\ve \to u
$$
in different topologies according to two different choices of the
kernel $G$.

\medskip

Let us give an heuristic idea in one space dimension, with $\Omega
= (0,1)$, of why the scaling involved in
\eqref{rescales} is the correct one. We assume that
$$\int^{\infty}_1 G(1,1-y)\,dy=-\int_{-\infty}^0 G(0,-y)\,dy
=\int_0^1J(y)\,y\,dy$$ and, as stated above, $G(x,\cdot)$ has
compact support independent of $x$.
In this case
\eqref{Neumann.C.1} reads
$$
\begin{array}{rl}
 u_t (x,t) =  & \displaystyle \frac{1}{\ve^2}\displaystyle \int_{0}^1
 J_\ve
 \left( x-y \right)\big(u(y,t) - u(x,t)\big)
  dy   + \displaystyle \frac{1}{\ve} \int_{-\infty}^0 G_\ve\left(x, x-y \right)
 g(y,t) \, dy \\[12pt]
& + \displaystyle \frac{1}{\ve} \int_{1}^{+\infty} G_\ve \left(x,
x-y\right)
 g(y,t) \, dy:= {\mathcal A}_\ep u (x,t).
 \end{array}
$$

If $x \in (0,1)$ a Taylor expansion gives that for any fixed
smooth $u$ and $\ep$ small enough, the right hand side  $
{\mathcal A}_\ep u$ in \eqref{Neumann.C.1} becomes
$$
{\mathcal A}_\ep u(x) = \frac{1}{\ve^2} \int_{0}^1 J_\ve \left(
x-y
\right) (u(y) - u(x) ) \, dy \approx u_{xx} (x)
$$
and if $x=0$ and $\ep$ small,
$$ {\mathcal A}_\ep u(0)=
\frac{1}{\ve^2} \int_{0}^1 J_\ve \left( -y \right) (u(y)
- u(0) ) \, dy  + \frac{1}{\ve} \int_{-\infty}^0 G_\ve
\left(0, -y \right) g(y)  \, dy \approx \frac{C_2}{\ve}(  u_x (0) - g(0)).
$$

Analogously, $ {\mathcal A}_\ep u(1)\approx (C_2 / \ve)(- u_x (1)
+ g(1))$. However, the proofs of our results are much more
involved than simple Taylor expansions due to the fact that for
each $\ep>0$ there are points $x\in \Omega$ for which the ball in
which integration takes place, $B(x,d\ve )$, is not contained in
$\Omega$. Moreover, when working in several space dimensions, one
has to take into account the geometry of the domain.

\medskip

Our first result deals with homogeneous boundary conditions, this
is, $g\equiv 0$.
\begin{teo}\label{zero} Assume $g\equiv 0$.
Let $\Omega$ be a bounded $C^{2+\alpha}$ domain for some
$0<\alpha<1$. Let $u\in C^{2+\alpha, 1+\alpha/2}(\overline
\Omega\times[0,T])$ be the solution to \eqref{calor.I}
and let $\uep$ be the solution to \eqref{Neumann.C.1} with $J_\ep$
as above. Then,
$$\sup_{t\in [0,T]}\|\uep (\cdot ,t)- u (\cdot
,t)\|_{L^{\infty}(\Omega )} \to 0$$ as $\ve \to 0$.
\end{teo}

Note that this result holds for every $G$ since $g\equiv 0$, and
that the assumed regularity in $u$ is standard if $u_0\in
C^{2+\alpha}(\overline\Omega)$ and $\partial u_0 /
\partial\eta=0$. See, for instance,
\cite{Friedman}.

We will prove Theorem \ref{zero} by constructing adequate super
and subsolutions and then using comparison arguments to get bounds
for the difference $u^\ve - u$.

Now we will make explicit the functions $G$ we will deal with in
order to consider $g\neq 0$.

To define the first one let us introduce some notation. As before,
let $\Omega$ be a bounded $C^{2+\alpha}$ domain. For $x
\in
\Omega_\ve := \{ x
\in
\Omega \ | \ \mbox{dist} (x,
\partial \Omega) < d \ve \}$ and $\ve$ small enough we write $x =
\bar{x} - s\, d \, \eta ( \bar{x})$ where $\bar{x}$ is the orthogonal
projection of $x$ on $\partial \Omega$, $0<s<\ve$ and $\eta
(\bar{x})$ is the unit exterior normal to $\Omega$ at $\bar x$.
Under these assumptions  we define
\begin{equation}\label{italiano}
G_1(x,\xi ) = -J(\xi)\, \eta (\bar{x}) \cdot \xi \quad\mbox{for }
x\in
\Omega_\ve.
\end{equation}

Notice that the last integral in \eqref{Neumann.C.1} only involves
points $x \in \Omega_\ve$ since when $y \not\in \Omega$, $x-y\in
supp\, J_\ep$ implies that $x\in\Omega_\ep$. Hence the above
definition makes sense for $\ep$ small.

For this choice of the kernel, $G=G_1$, we have the following
result.

\begin{teo}\label{conv} Let $\Omega$ be a bounded $C^{2+\alpha}$
domain, $g \in C^{1+\alpha,
(1+\alpha)/2}(\partial\Omega\times[0,T])$, $u\in C^{2+\alpha,
1+\alpha/2}(\overline\Omega\times[0,T])$ the solution to
\eqref{calor.I}, for
some $0<\alpha<1$. Let $J$ as before and
 $G(x,\xi )=G_1(x,\xi )$, where $G_1$ is defined by \eqref{italiano}. Let $\uep$ be the solution to
 \eqref{Neumann.C.1}. Then,
 $$\sup_{t\in [0,T]}\|\uep (\cdot
,t)- u (\cdot ,t)\|_{L^1(\Omega )} \to 0$$ as $\ve \to 0$.
\end{teo}

Observe that $G_1$ may fail to be nonnegative and hence a
comparison principle may not hold. However, in this case our proof
of convergence to the solution of the heat equation does not rely
on comparison arguments for \eqref{Neumann}. If we want a
nonnegative kernel $G$, in order to have a comparison principle,
we can modify $(G_1)_\ep$ by taking instead
$$
(\tilde{G}_1)_\ep (x, \xi) = (G_1)_\ve (x, \xi) +\kappa
\varepsilon J_\varepsilon (\xi) = \frac{1}{\ve} J_\ve (\xi)
\left(- \eta
(\bar x) \cdot \xi +\kappa \ve^{2}
\right).
$$
Note that for $x\in\overline\Omega$ and $y\in\Oc$,
$(\tilde{G}_1)_\ep (x, x-y)=\frac{1}{\ve} J_\ve (x-y) \left(- \eta
(\bar x) \cdot (x-y) +\kappa \ve^{2}
\right)$  is nonnegative for $\varepsilon$ small if we choose the constant $\kappa$ as
a bound for the curvature of $\partial \Omega$, since $|x-y| \le d
\, \ve$. As will be seen in Remark~\ref{suspiro.lima},
Theorem~\ref{conv} remains valid with $(G_1)_\ep$ replaced by
$(\tilde{G}_1)_\ep$.

Finally, the other ``Neumann'' kernel we propose is
$$
G(x,
\xi) = G_2 (x, \xi) =  C_2 J(\xi),
$$
where $C_2$ is such that
\begin{equation}\label{C2}
\int_0^d\int_{\{z_N>s\}}J(z)\big(C_2 - z_N\big)\,dz\,ds=0.
\end{equation}
This choice of $G$ is natural since we are considering a flux with
a jumping probability that is a scalar multiple of the same
jumping probability that moves things in the interior of the
domain, $J$.

Several properties of solutions to
\eqref{Neumann} have been recently investigated in \cite{cerw} in
the case $G=G_2$ for different choices of $g$.

For the case of $G_2$ we can still prove convergence but in a
weaker sense.

\begin{teo}\label{teo.conv.J} Let $\Omega$ be a bounded $C^{2+\alpha}$
domain, $g \in C^{1+\alpha,
(1+\alpha)/2}(\partial\Omega\times[0,T])$, $u\in C^{2+\alpha,
1+\alpha/2}(\overline\Omega\times[0,T])$ the solution to
\eqref{calor.I}, for
some $0<\alpha<1$. Let $J$ as before and
 $G(x,\xi )=G_2(x,\xi )= C_2 J(\xi )$, where $C_2$ is defined by \eqref{C2}. Let $\uep$ be the solution to
 \eqref{Neumann.C.1}. Then,
 for each $t\in [0,T]$
$$ u_\ep(x,t)\rightharpoonup u(x,t)\quad*-\mbox{weakly in }
L^\infty(\Omega) $$ as $\ve \to 0$.
\end{teo}

The rest of the paper is organized as follows: in
Section~\ref{Sect.exist.uni} we prove existence, uniqueness and a
comparison principle for our nonlocal equation. In
Section~\ref{Sect.Proof.teo.1} we prove the uniform convergence
when $g=0$. In Section~\ref{Sect.L1} we deal with the case $G=G_1$
and finally in Section~\ref{sect.weak.L1} we prove our result when
$G=G_2$.

\section{Existence and uniqueness}
\label{Sect.exist.uni} \setcounter{equation}{0}

In this section we deal with existence and uniqueness of solutions
of \eqref{Neumann}. Our result is valid in a general $L^1$
setting.

\begin{teo} \label{teo.Neumann.Lineal} Let $\Omega$ be a bounded domain.
Let $J\in L^1(\RR^N)$ and  $G\in L^\infty(\Omega\times\RR^N)$. For
every $u_0 \in L^1(\Omega)$ and $g\in L^\infty_{loc}
([0,\infty);L^1(\Oc))$ there exists a unique solution $u$ of
\eqref{Neumann} such that $u \in C ([0, \infty ); L^1(\Omega))$
and $u(x,0)=u_0(x)$.
\end{teo}

As in \cite{cer} and \cite{cerw}, existence and uniqueness will be
a consequence of Banach's fixed point theorem. We follow closely
the ideas of those works in our proof, so we will only outline the
main arguments. Fix $t_0
>0$ and consider the Banach space
$$X_{t_0} =
C ([0,t_0]; L^1(\Omega))$$ with the norm
$$ ||| w|||=  \max\limits_{0\le t \le t_0} \| w(\cdot,
t) \|_{L^1(\Omega)}. $$ We will obtain the solution as a fixed
point of the operator $T_{u_0,g}: X_{t_0} \to X_{t_0}$ defined by
\begin{equation}\label{T}
\begin{array}{rl}
\displaystyle T_{u_0,g} (w) (x,t) =  u_0 (x) & + \displaystyle
\int_0^t \int_\Omega J\left( x-y \right) (w(y,s) - w(x,s))\, dy \,
ds \\[12pt]
& + \displaystyle \int_0^t \int_{\RR^N\setminus\Omega} G(x,x-y)
 g(y,t) \, dy \, ds.
\end{array}
\end{equation}

The following lemma is the main ingredient in the proof of
existence.

\begin{lema}\label{contraction} Let $J$ and $G$ as in Theorem \ref{teo.Neumann.Lineal}.
Let $g,\ h \in L^\infty ((0,t_0);L^1(\Oc))$ and $u_0,\ v_0 \in L^1
(\Omega)$. There exists a constant $C$ depending only on $\Omega$,
$J$ and $G$ such that for $w,z \in X_{t_0}$,
\begin{equation}
\label{ec.contrac}
|||T_{u_0,g} (w)-T_{v_0,h} (z)|||
\leq \|u_0-v_0\|_{L^1}+ C t_0 \left( |||w-z||| +\|g-h\|_{L^\infty((0,t_0);L^1(\Oc))}\right).
\end{equation}
\end{lema}

\begin{proof} We have
$$\begin{aligned} \displaystyle
& \int_\Omega |T_{u_0,g} (w)(x,t) - T_{v_0,h}(z)(x,t)|\, dx
\le\int_\Omega|u_0(x)-v_0(x)|\,dx
 \\
& \displaystyle \qquad +  \int_\Omega \left| \int_0^t
\int_\Omega J\left( x-y \right) \Big[ (w(y,s) - z(y,s))  - (w(x,s)
- z(x,s)) \Big]
 \, dy \, ds \right| \, dx\\[8pt]
 &\displaystyle \qquad +\int_\Omega \int_0^t
\int_{\Oc} |G(x,x-y)||g(y,s)-h(y,s)|\,dy\,ds\,dx.
\end{aligned}$$
Therefore, we obtain \eqref{ec.contrac}.
\end{proof}

\begin{proof}[Proof of Theorem \ref{teo.Neumann.Lineal}] Let $T=T_{u_0,g}$.
We check first that $T$ maps $X_{t_0}$ into $X_{t_0}$. From
\eqref{T} we see that for $0\le t_1<t_2\le t_0$,
$$
\|T(w)(t_2)-T(w)(t_1)\|_{L^1(\Omega)}\le
A\int_{t_1}^{t_2} \int_\Omega  |w(y,s)| \, dy \,
ds+B\int_{t_1}^{t_2}\int_{\RR^N\setminus\Omega} |g(y,s)|  \, dy
\, ds.
$$
On the other hand, again from \eqref{T}
$$
\|T(w)(t)-u_0\|_{L^1(\Omega)}\le C t\big\{|||w||| +\|g\|_{
L^\infty ((0,t_0);L^1(\Oc))}\big\}.
$$
These two estimates give that $T(w)\in C([0,t_0]; L^1(\Omega))$.
Hence $T$ maps $X_{t_0}$ into $X_{t_0}$.

Choose $t_0$ such that $ C t_0 <1$. {}From Lemma \ref{contraction}
we get that $T$ is a strict contraction in $X_{t_0}$ and the
existence and uniqueness part of the theorem follows from Banach's
fixed point theorem in the interval $[0,t_0]$. To extend the
solution to $[0,\infty)$ we may take as initial datum $u(x,t_0)\in
L^1(\Omega)$ and obtain a solution in $[0,2\, t_0]$. Iterating
this procedure we get a solution defined in $[0,\infty)$.
\end{proof}

Our next aim is to prove a comparison principle for
\eqref{Neumann} when $J,\  G\ge 0$. To this end we define what we understand by
sub and supersolutions.

\begin{defi} A function $u\in C([0,T);L^1((\Omega))$ is a
supersolution of \eqref{Neumann} if $u(x,0) \ge u_0(x)$ and
$$
u_t(x,t)\ge \int_\Omega J(x-y)\big(u(y,t)-u(x,t)\big)\,dy+
\int_{\Oc} G(x, x-y) g(y,t)\,dy.
$$
\end{defi}

Subsolutions are defined analogously by reversing the
inequalities.

\begin{lema} \label{compar.0} Let $J,\ G\ge 0$, $u_0\ge 0$ and $g\ge0$. If $u\in
C(\overline{\Omega}\times [0,T])$ is a supersolution to
\eqref{Neumann}, then $u\ge0$.
\end{lema}

\begin{proof}
Assume that $u(x,t)$ is negative somewhere. Let $v(x,t)=u(x,t)+\ep
t$ with $\ep$ so small such that $v$ is still negative somewhere.
Then, if we take $(x_0,t_0)$ a point where $v$ attains its
negative minimum, there holds that $t_0>0$ and
$$
\begin{aligned}
&v_t(x_0,t_0)=u_t(x_0,t_0)+\ep>\int_{\Omega} J(x-y)(u(y,t_0) -
u(x_0,t_0)) \, dy\\
&\hskip1cm=\int_{\Omega} J(x-y)(v(y,t_0) - v(x_0,t_0)) \, dy\ge0
\end{aligned}
$$
which is a contradiction. Thus, $u\ge0$.
\end{proof}

\begin{coro} \label{coro.compar} Let $J,\ G \ge 0$ and bounded. Let $u_0$
and $v_0$ in $L^1(\Omega)$ with $u_0\ge v_0$ and $g,\,h\in
L^\infty((0,T);L^1(\RR^N\setminus\Omega))$ with $g\ge h$. Let $u$
be a solution of \eqref{Neumann} with initial condition $u_0$ and
flux $g$ and $v$ be a solution of
\eqref{Neumann} with initial condition $v_0$ and flux $h$. Then,
$$u\ge v \qquad \mbox{ a.e}.$$
\end{coro}

\begin{proof}
Let $w=u-v$. Then, $w$ is a supersolution with initial datum
$u_0-v_0 \ge 0$ and boundary datum $g-h\ge 0$. Using the
continuity of solutions with respect to the initial and Neumann
data (Lemma \ref{contraction}) and the fact that $J\in
L^\infty(\RR^N)$, $G\in L^\infty(\Omega\times\RR^N)$ we may assume
that $u,\, v \in C(\overline{\Omega}\times [0,T])$. By Lemma
\ref{compar.0} we obtain that $w=u-v \ge 0$. So the corollary is
proved.
\end{proof}

\begin{coro}  Let $J,\ G \ge 0$ and bounded. Let  $u\in C(\overline{\Omega}\times [0,T])$ (resp. $v$) be a
supersolution (resp. subsolution) of \eqref{Neumann}. Then, $u\ge
v$.
\end{coro}

\begin{proof} It follows the lines of the proof of the previous corollary.
\end{proof}

\section{Uniform convergence in the case $g\equiv 0$}
\label{Sect.Proof.teo.1} \setcounter{equation}{0}

In order to prove  Theorem \ref{zero} we set $w^\ve = \uep -u$ and
let $\tilde{u}$ be a $C^{2+\alpha, 1+\alpha/2}$ extension of $u$
to $\mathbb{R}^N\times[0,T]$. We define
$$L_\ve (v) = \frac{1}{\ve^2}\int_{\Omega}J_\ve (x-y)\big(v(y,t)-v(x,t)\big)dy$$
and
$$\tilde{L}_\ve (v) = \frac{1}{\ve^2}\int_{\mathbb{R}^N}J_\ve
(x-y)\big(v(y,t)-v(x,t)\big)dy.$$

Then
$$
\begin{aligned}
\wep _t &=L_\ve (\uep) - \Delta u +
\frac{1}{\ve}\int_{\Oc}G_\ve (x,x-y)g(y,t)\, dy  \\
&=L_\ve (\wep) +\tilde{L}_\ve (\tilde u) - \Delta u
+\frac{1}{\ve}\int_{\Oc}G_\ve (x,x-y)g (y,t)\, dy  \\
&\hskip2cm-\frac{1}{\ve^2}\int_{\Oc}J_\ve
(x-y)\big(\tilde{u}(y,t)-\tilde{u}(x,t)\big)\, dy.
\end{aligned}
$$
Or
$$
\wep_t-L_\ve (\wep) = F_\ve (x,t),
$$
where, noting that $\Delta u = \Delta \tilde u$ in $\Omega$,
$$
\begin{aligned}
F_\ve (x,t) =  \displaystyle
\tilde{L}_\ve (\tilde u) - \Delta \tilde u
&+\frac{1}{\ve}\int_{\Oc}G_\ve (x,x-y)g (y,t)\, dy \\
&\displaystyle -\frac{1}{\ve^2}\int_{\Oc}J_\ve
(x-y)\big(\tilde{u}(y,t)-\tilde{u}(x,t)\big)\, dy.
\end{aligned}
$$

Our main task in order to prove the uniform convergence result is
to get bounds on $F_\ve$.

First, we observe that it is well known that by the choice of
$C_1$, the fact that $J$ is radially symmetric and $\tilde u \in
C^{2+\alpha, 1+\alpha/2}(\mathbb{R}^N\times[0,T])$, we have that
\begin{equation}\label{est1}
\sup_{t\in [0,T]}\|\tilde{L}_\ve (\tilde u) - \Delta \tilde u
\|_{L^\infty (\Omega) } = O(\ve^{\alpha} )\, .
\end{equation}

In fact,
$$
\frac{C_1}{\ve^{N+2}} \int_{\RR^N} J \left( \frac{x-y}{\ve} \right)
\left(\tilde{u} (y,t) - \tilde{u} (x,t)\right)\, dy - \Delta \tilde{u} (x,t)
$$
 becomes, under the change variables  $z= (x-y) / \ve$,
$$
\frac{C_1}{\ve^{2}} \int_{\RR^N} J \left( z \right)
(\tilde{u} (x-\ve z,t) - \tilde{u} (x,t))\, dy - \Delta \tilde{u}
(x,t)
$$
and hence \eqref{est1} follows by a simple Taylor expansion.

\medskip

Next, we will estimate the last integral in $F_\ep$. We remark
that the next lemma is valid for any smooth function, not only for
a solution to the heat equation.

\begin{lema}\label{F}
If $\theta$  is a $C^{2+\alpha, 1+\alpha/2}$ function  on
$\mathbb{R}^N\times[0,T]$ and $\di\frac{\partial
\theta}{\partial\eta}= h $ on $\partial\Omega$,  then for
$x\in\Omega_\ep=
\{z\in\Omega\,|\,{\rm dist\,}(z,\partial\Omega)<d\ep\}$ and $\ep$ small,
$$
\begin{aligned}
 &\displaystyle \frac{1}{\ve^2}\int_{\Oc}J_\ve
(x-y)\big(\theta (y,t)-\theta (x,t)\big)\, dy = \frac{1}{\ve}
\int_{\Oc} J_\ve (x-y) \eta(\bar x)\cdot
\frac{(y-x)}{\ve} h (\bar x,t)\,dy
\\
& \qquad +\int_{\Oc}J_\ve (x-y)\sum_{|\beta|=2}
\frac{D^{\beta}\theta }{2}(\bar x,t)\Big[\big(\frac{(y-\bar
x)}{\ep}\big)^{\beta}-\big(\frac{( x-\bar
x)}{\ep}\big)^{\beta}\Big]dy
 +O(\ve^{\alpha} ),
\end{aligned}
$$
where $\bar x$ is the orthogonal projection of $x$ on the boundary
of $\Omega$ so that $\|\bar x-y\|\le 2d\ep$.
\end{lema}

\begin{proof}
 Since $\theta \in C^{2+\alpha,
1+\alpha/2}(\RR^N \times[0,T] )$ we have
$$
\begin{aligned}
\theta (y,t)-\theta (x,t)& = \theta (y,t)-\theta (\bar x,t)
-\big(\theta (x,t)-\theta (\bar x,t)\big)\\
&=\nabla \theta (\bar x,t)
\cdot (y- x)
+ \sum_{|\beta|=2}\frac {D^{\beta}\theta}{2} (\bar
x,t)\big[(y-\bar x)^{\beta}-(x-\bar x)^{\beta}\big]\\
&\qquad +O(||\bar x-x||^{2+\alpha}) + O(||\bar x-y||^{2+\alpha}).
\end{aligned}
$$

Therefore,
$$
\begin{aligned}
& \frac{1}{\ve^2}\int_{\Oc}J_\ve (x-y)\big(\theta (y,t)-\theta
(x,t)\big)\, dy =\frac{1}{\ve}\int_{\Oc}J_\ve (x-y)\nabla
\theta (\bar x,t)
\cdot \frac{(y-x)}{\ve} \, dy\\
& \quad +\int_{\Oc}J_\ve
(x-y)\sum_{|\beta|=2}\frac{D^{\beta}\theta}{2}(\bar
x,t)\Big[\big(\frac{(y-\bar x)}{\ep}\big)^{\beta} -\big(\frac{(
x-\bar x)}{\ep}\big)^{\beta}\Big]\, dy+ O(\ve^{\alpha} ).
\end{aligned}
$$

Fix $x\in \Omega_\ep$. Let us take a new coordinate system such
that $\eta(\bar x)=e_N$. Since
$\di\frac{\partial\theta}{\partial\eta}=h $ on $\partial\Omega$,
we get
$$
\begin{aligned}
&  \int_{\Oc} J_\ve (x-y)\nabla \theta (\bar x,t) \cdot
\frac{(y-x)}{\ve} \, dy\\
 & = \int_{\Oc} J_\ve (x-y) \eta(\bar x)\cdot \frac{(y-x)}{\ve} h
(\bar x,t) \, dy + \int_{\Oc}J_\ve (x-y)\sum_{i=1}^{N-1}
\theta_{x_i}(\bar x,t) \di\frac{(y_i-x_i)}{\ve}\, dy.
\end{aligned}
$$

We will estimate this last integral. Since $\Omega$ is a
$C^{2+\alpha}$ domain we can chose vectors $e_1$, $e_2$, ...,
$e_{N-1}$ so that there exists $\kappa>0$ and constants
$f_{i}(\bar x)$ such that
$$\begin{aligned} &B_{2d\ep}(\bar x)\cap\left\{y_N-\big(\bar x_N
+\sum_{i=1}^{N-1}f_{i}(\bar
x)(y_i- x_i)^2\big)> \kappa \ep^{2+\alpha} \right\} \subset \Oc,\\
&B_{2d\ep}(\bar x)\cap\left\{y_N-\big(\bar x_N
+\sum_{i=1}^{N-1}f_{i}(\bar x)(y_i- x_i)^2\big)<-\kappa
\ep^{2+\alpha} \right\} \subset \Omega.
\end{aligned}
$$
Therefore
$$
\begin{aligned}
&\int_{\Oc}J_\ve (x-y)\Big(
 \sum_{i=1}^{N-1} \theta_{x_i}(\bar
x,t)
\di\frac{(y_i-x_i)}{\ve}\Big)\, dy\\
&\hskip1cm=\int_{(\Oc)\cap\left\{\big|y_N-\big(\bar x_N
+\sum_{i=1}^{N-1}f_{i}(\bar x)(y_i- x_i)^2\big)\big|\leq \kappa
\ep^{2+\alpha} \right\}}J_\ve (x-y)\Big( \sum_{i=1}^{N-1}
\theta_{x_i}(\bar x,t)
\di\frac{(y_i-x_i)}{\ve}\Big)\, dy\\
&\hskip1cm \qquad +\int_{\left\{y_N-\big(\bar x_N
+\sum_{i=1}^{N-1}f_{i}(\bar x)(y_i- x_i)^2\big)>
\kappa\ep^{2+\alpha} \right\}}J_\ve (x-y)\Big(
\sum_{i=1}^{N-1} \theta_{x_i}(\bar x,t)
\di\frac{(y_i-x_i)}{\ve}\Big)\, dy\\
&\hskip1cm=I_1+I_2.
\end{aligned}
$$
If  we take $z=(y-x)/\ep$ as a new variable, recalling that
$\bar{x}_N- x_N = \ep s$, we obtain
$$
|I_1|\le C_1\sum_{i=1}^{N-1} |\theta_{x_i}(\bar x,t)|
\int_{\left\{\big|z_N-\big (s+\ep\sum_{i=1}^{N-1}f_{i}(\bar
x)(z_i)^2\big)\big|\le \kappa\ep^{1+\alpha} \right\}} J(z) |z_i|
\, dz\le C\, \kappa\, \ep^{1+\alpha}.
$$
On the other hand,
$$
\begin{aligned}
I_2 &= C_1\sum_{i=1}^{N-1} \theta_{x_i}(\bar x,t)
\int_{\left\{z_N-\big (s+\ep\sum_{i=1}^{N-1}f_{i}(\bar
x)(z_i)^2\big)> \kappa\ep^{1+\alpha} \right\}} J(z)\, z_i \,  dz.\\
\end{aligned}
$$
Fix $1\le i\le N-1$. Then, since $J$ is radially symmetric, $J(z)
\,z_i$ is an odd function of the variable $z_i$ and, since
the set $\left\{z_N-\big (s+\ep\sum_{i=1}^{N-1}f_{i}(\bar
x)(z_i)^2\big)> \kappa\ep^{1+\alpha} \right\}$  is symmetric in
that variable we get
$$
I_2=0.
$$
Collecting the previous estimates the lemma is proved.
\end{proof}

We will also need the following inequality.

\begin{lema}\label{lomitos}
There exist $K>0$ and $\bar{\ve}>0$ such that, for $\ep
<\bar{\ve}$,
\begin{equation}\label{integrales.II}
\int_{\Oc} J_\ep(x-y)\eta(\bar x)\cdot \frac{(y-x)}{\ve}\,dy \ge K \int_{\Oc}J_\ve (x-y)\,dy.
\end{equation}
\end{lema}

\begin{proof}
Let us put the origin at the point $\bar x$  and take a coordinate
system such that $\eta(\bar x)= e_N$. Then, $x=(0,-\mu)$ with
$0<\mu<d\ep$. Then, arguing as before,
$$
\begin{aligned}
&\int_{\Oc} J_\ep(x-y)\eta(\bar x)\cdot \frac{(y-x)}{\ve}\,dy=\int_{\Oc} J_\ep(x-y)\frac{y_N+\mu}\ep\,dy\\
&= \int_{\{y_N>\kappa
\ep^2\} }J_\ep(x-y)\frac{y_N+\mu}\ep\,dy+\int_{\Oc\cap\{|y_N|<\kappa
\ep^2\}}
J_\ep(x-y)\frac{y_N+\mu}\ep\,dy\\
&\ge \int_{\{y_N>\kappa
\ep^2\}}J_\ep(x-y)\frac{y_N+\mu}\ep\,dy-C\ep.
\end{aligned}
$$

Fix $c_1$ small such that
$$
\frac{1}{2}\int_{\{z_N>0\}}J(z) \,z_N\,dz \ge 2
c_1\int_{\{0<z_N<2c_1\}}J(z)\,dz.
$$
We divide our arguments into two cases according to whether
$\mu\le c_1 \ep$ or $\mu > c_1 \ep$.

\noindent {\bf Case I} Assume $\mu \le c_1 \ep$. In this case we have,
\begin{equation}\label{c1}
\begin{aligned}
&\int_{\{y_N>\kappa \ep^2\}}J_\ep(x-y)\frac{y_N+\mu}\ep\,dy=
C_1\int_{\{z_N>\kappa \ep+\frac\mu\ep\}}J(z) \,z_N\,dz\\
&= C_1\left( \int_{\{z_N>0\}}J(z) \,z_N\,dz-\int_{\{0<z_N<\kappa \ep+\frac\mu\ep\}}J(z) \,z_N\,dz \right)\\
& \ge C_1\left( \int_{\{z_N>0\}}J(z) \,z_N\,dz-2
c_1\int_{\{0<z_N<2c_1\}}J(z)\,dz\right)\ge\frac{C_1}2
\int_{\{z_N>0\}}J(z)
\,z_N\,dz.
\end{aligned}
\end{equation}

Then,
$$
\begin{aligned}
& \int_{\Oc} J_\ep(x-y)\eta(\bar x)\cdot \frac{(y-x)}{\ve}\,dy- K \int_{\Oc}J_\ve (x-y)\,dy \\
&\ge C_1\left( \frac{1}{2 }\int_{\{z_N>0\}}J(z) \,z_N\,dz-
K\right)-C\ve \ge 0,
\end{aligned}
$$
if $\ep$ is small enough and
$$
K< \frac{1}{4 }\int_{\{z_N>0\}}J(z) \,z_N\,dz .
$$

\medskip

\noindent {\bf Case II} Assume that $\mu\ge c_1 \ep$. For $y$
in $\RR^N \setminus \Omega \cap B(\bar x, d\ve)$ we have
$$
\frac{y_N}{\varepsilon} \ge -\kappa \varepsilon.
$$
Then,
$$
\begin{aligned}
& \int_{\Oc} J_\ep(x-y)\frac{y_N+\mu}\ep\,dy- K \int_{\Oc}J_\ve (x-y)\,dy \\
& \ge  (c_1-\kappa\ep)\int_{\Oc} J_\ep(x-y)\,dy- K \int_{\Oc}J_\ve (x-y)\,dy \\
& = \big( c_1 - \kappa\ep -K \big)\int_{\Oc}J_\ve (x-y)\,dy
\ge 0,
\end{aligned}
$$
if $\ep$ is small and
$$
K < \frac{c_1}{2}.
$$
This ends the proof of \eqref{integrales.II}.
\end{proof}

We now prove Theorem \ref{zero}.

\begin{proof}[Proof of Theorem \ref{zero}] We will use a comparison argument. First, let us
look for a supersolution. Let us pick an auxiliary function $v$ as
a solution to
$$
\begin{cases}
v_t-\Delta v=h(x,t)\quad&\mbox{in}\quad \Omega\times(0,T),\\[10pt]
\di\frac{\partial v}{\partial\eta}=g_1(x,t)\quad&\mbox{on}\quad
\partial\Omega\times(0,T),\\[10pt]
v(x,0)=v_1(x)\quad&\mbox{in}\quad \Omega.
\end{cases}
$$
for some smooth functions $h(x,t)\geq 1$, $g_1(x,t)\geq 1$ and
$v_1(x)\geq 0$ such that the resulting $v$ has an extension
$\tilde v$ belongs to $C^{2+\alpha, 1+\alpha/2}(\mathbb{R}^N\times
[0,T])$, and  let $M$ be an upper bound for $v$ in $\bar
\Omega\times [0,T]$. Then,
$$
v_t= L_\ep v +(\Delta v- \tilde{L}_\ve  \tilde v )+
\frac{1}{\ve^2}\int_{\Oc}J_\ve
(x-y)(\tilde{v}(y,t)-\tilde{v}(x,t))\,dy+ h(x,t).
$$

Since $\Delta v = \Delta \tilde v$ in $\Omega$, we have that $v$
is a solution to
$$
\begin{cases}
v_t-L_\ep v=H(x,t,\ep)\quad&\mbox{in}\quad \Omega\times(0,T),\\
v(x,0)= v_1(x)\quad&\mbox{in}\quad \Omega,
\end{cases}
$$
where by \eqref{est1}, Lemma \ref{F}  and the fact that $h\ge 1$,
$$
\begin{aligned}
H(x,t,\ep)=& (\Delta
\tilde v- \tilde{L}_\ve \tilde v )+ \frac{1}{\ve^2}\int_{\Oc}J_\ve
(x-y)(\tilde{v}(y,t)-\tilde{v}(x,t))\,dy+ h(x,t)\\
\geq &
\Big(\frac{1}{\ve} \int_{\Oc} J_\ve
(x-y)\eta(\bar x)\cdot \frac{(y-x)}{\ve} g_1(\bar
x,t)\,dy\\
& +\int_{\Oc}J_\ve
(x-y)\sum_{|\beta|=2}\frac{D^{\beta}\tilde{v}}{2}(\bar
x,t)\Big[\big(\frac{(y-\bar x)}{\ep}\big)^{\beta}
-\big(\frac{( x-\bar x)}{\ep}\big)^{\beta}\Big]\,dy \Big)+ 1-C \ep^{\alpha}\\
\geq & \Big(\frac{ g_1(\bar
x,t)}{\ve} \int_{\Oc} J_\ve (x-y) \eta(\bar x)\cdot
\frac{(y-x)}{\ve}\,dy -D_1\int_{\Oc}J_\ve (x-y)\, dy
\Big)+\frac{1}{2}
\end{aligned}
$$
for some constant $D_1$ if $\ep$ is small so that $C\ep^\alpha\le
1/ 2$.

Now, observe that Lemma \ref{lomitos} implies that for every
constant $C_0>0$ there exists $\ep_0$ such that,
$$
\frac{1}{\ep}\int_{\Oc} J_\ep(x-y)\eta(\bar x)\cdot
\frac{(y-x)}{\ve}\,dy - C_0 \int_{\Oc}J_\ve (x-y)\,dy \ge 0,
$$
if $\ep<\ep_0$.

Now, since $g=0$, by \eqref{est1} and Lemma \ref{F} we obtain
$$
\begin{aligned}
|F_\ep|&\le C\ep^{\alpha}+ \int_{\Oc}J_\ve
(x-y)\sum_{|\beta|=2}\frac{D^{\beta}\tilde{u}}{2}(\bar
x,t)\Big[\big(\frac{(y-\bar x)}{\ep}\big)^{\beta}
-\big(\frac{( x-\bar x)}{\ep}\big)^{\beta}\Big]dy\\
& \leq  C\ep^{\alpha}+ C_2\int_{\Oc}J_\ve (x-y)\,dy.
\end{aligned}
$$

Given $\delta >0$, let $v_\delta = \delta v$. Then $v_\delta$
verifies
$$
\begin{cases}
(v_\delta)_t-L_\ep v_\delta=\delta H(x,t,\ep)\quad&\mbox{in}\quad \Omega\times(0,T),\\
v_\delta(x,0)= \delta v_1(x)\quad&\mbox{in}\quad \Omega.
\end{cases}
$$

By our previous estimates, there exists $\ep_0=\ep_0(\delta)$ such
that for $\ep\le\ep_0$,
$$
|F_\ep|\le \delta H(x,t,\ep).
$$
So, by the comparison principle for any $\ep
\leq \ep_0$ it holds that
$$-M\delta \leq -v_\delta \leq w_\ep \leq v_\delta\leq M \delta.$$

Therefore, for every $\delta>0$,
$$-M\delta \leq \liminf_{\ep\to0}w_\ep \leq\limsup _{\ep\to0}w_\ep \leq M \delta.$$
and the theorem is proved.
\end{proof}

%%%%%%%%%%%%%%%%%%%%%%%%%%%%%%%%%%%%%%%%%%%%%%%%%%%%%%%%%%%%%
%%%%%%%%%%%%%%%%%%%%%%%%%%%%%%%%%%%%%%%%%%%%%%%%%%%%%%%%%%%%%%
%%%%%%%%%%%%%%%%%%%%%%%%%%%%%%%%%%%%%%%%%%%%%%%%%%%%%%%%%%%%%%%
%%%%%%%%%%%%%%%%%%%%%%%%%%%%%%%%%%%%%%%%%%%%%%%%%%%%%%%%%%%%%%%%

\section{Convergence in $L^1$ in the case $G=G_1$} \label{Sect.L1}
\setcounter{equation}{0}

First we prove that $F_\ve$ goes to zero as $\ve $ goes to zero.

\begin{lema}\label{fep}
If $G=G_1$ then
 $$ F_\ve (x,t) \rightarrow 0\quad\mbox{ in}\quad
L^{\infty}\big([0,T];L^1(\Omega) \big) $$ as $\ve \to 0$.
\end{lema}

\begin{proof}
As $G=G_1=  -J(\xi)\, \eta (\bar{x}) \cdot \xi$, for  $x \in
\Omega_\ve$, by \eqref{est1} and Lemma~\ref{F},
$$
\begin{aligned}
F_\ve (x,t) = & \frac{1}{\ep}  \int_{\Oc}J_\ve (x-y) \eta(\bar
x)\cdot
\frac{(y-x)}{\ve} \big(g(y,t)-  g(\bar x,t)\big)\,dy\\
& - \int_{\Oc}J_\ve
(x-y)\sum_{|\beta|=2}\frac{D^{\beta}\tilde{u}}{2}(\bar
x,t)\Big[\big(\frac{(y-\bar x)}{\ep}\big)^{\beta} -\big(\frac{(
x-\bar x)}{\ep}\big)^{\beta}\Big]\,dy+ O(\ve^{\alpha} ).
\end{aligned}
$$
As $g$ is smooth, we have that $F_\ep$ is bounded in $\Omega_\ve$.
Recalling the fact that $|\Omega_\ve | = O(\ve )$ and $F_\ve
(x,t)=O(\ve^{\alpha} ) $ on $\Omega\setminus\Omega_\ve$ we get the
convergence result.
\end{proof}

We are now ready to prove Theorem \ref{conv}.

\begin{proof}[Proof of Theorem \ref{conv}] In the case $G=G_1$
we have proven in Lemma \ref{fep} that $F_\ep\to 0$ in
$L^1(\Omega\times[0,T])$. On the other hand, we have that $w^\ep =
u^\ep - u$ is a solution to
$$
\begin{aligned}
w_t-L_\ep (w)&=F_\ep\\
w(x,0)&=0.
\end{aligned}
$$
Let ${z}^\ep$ be a solution to
$$
\begin{aligned}
z_t-L_\ep (z)&=|F_\ep |\\
z(x,0)&=0.
\end{aligned}
$$
Then $-{z}^\ep$ is a solution to
$$
\begin{aligned}
z_t-L_\ep (z)&=-|F_\ep |\\
z(x,0)&=0.
\end{aligned}
$$
By comparison we have that
$$
-{z}^\ep \leq w^\ep \leq {z}^\ep \quad\mbox{and}\quad  {z}^\ep
\geq 0.
$$

Integrating  the equation for $ z^\ep$ we get
$$\|{z}^\ve (\cdot ,t)\|_{L^1(\Omega )} =     \int_{\Omega} {z}^\ve
(x,t)\, dx =
\int_{\Omega}\int_0^t|F_\ep (x,s)|\;ds\,dx\, .$$
Applying Lemma \ref{fep} we get
$$\sup_{t\in [0,T]}\|{z}^\ve(\cdot ,t) \|_{L^1(\Omega  )} \to 0$$
as $\ve \to 0$. So the theorem is proved.
\end{proof}

\begin{rem}\label{suspiro.lima} Notice that if we consider a
kernel $G$ which is a modification of $G_1$ of the form
$$
G_\ve (x,\xi) = (G_1)_\ve (x, \xi) + A (x, \xi , \ve)
$$
with
$$
\int_{\RR^N \setminus \Omega} |A(x, x-y ,\ve)| \, dy\to 0
$$
in $L^1 (\Omega)$ as $\ve \to 0$, then the conclusion of Theorem
\ref{conv} is still valid. In particular, we can take $A(x, \xi, \ve) =
\kappa\ve J_\ve (\xi)$.
\end{rem}

\section{Weak convergence in $L^1$ in the case $G=G_2$}
\label{sect.weak.L1} \setcounter{equation}{0}

First, we prove that in this case $F_\ve$ goes to zero as
measures.

\begin{lema} \label{lema.3.2.II}
If $G=G_2$ then there exists a constant $C$ independent of $\ve$
such that
$$
\int_0^T
\int_\Omega |F_\ve (x,s)| \, dx \, ds \le C.
$$
Moreover,
 $$
 F_\ve (x,t) \rightharpoonup 0\quad\mbox{ as measures}
 $$
 as $\ve \to 0$.
 That is, for any continuous function $\theta$, it holds that
$$
\int_0^T\int_{\Omega} F_\ep(x,t)\theta
(x,t)\,dx\,dt \to 0
$$
as $\ep \to 0$.
\end{lema}

\begin{proof}
As $G=G_2=C_2 J(\xi)$ and $g$ and $\tilde{u}$ are smooth, taking
again the coordinate system of Lemma
\ref{F}, we obtain
$$
\begin{aligned}
F_\ve (x,t) =&\frac{1}{\ep} \int_{\Oc}J_\ve (x-y)\Big (C_2g(y,t)-  \frac{y_N-x_N}{\ve} g(\bar x,t)\Big)\\
& - \frac{1}{\ep} \int_{\Oc} J_\ve (x-y) \sum_{i=1}^{N-1}
\tilde{u}_{x_i}(\bar x,t)
 \di\frac{(y_i-x_i)}{\ve}\, dy\\
&- \int_{\Oc}J_\ve
(x-y)\sum_{|\beta|=2}\frac{D^{\beta}\tilde{u}(\bar
x,t)}2\Big[\big(\frac{(y-\bar x)}{\ep}\big)^{\beta}
-\big(\frac{( x-\bar x)}{\ep}\big)^{\beta}\Big]dy+ O(\ve^{\alpha} )\\
= & \frac{1}{\ep}
\int_{\Oc}J_\ve (x-y)\Big(C_2g(\bar x,t) - \frac{y_N-x_N}{\ve} g(\bar x,t)\Big)\\
& - \frac{1}{\ep} \int_{\Oc}J_\ve (x-y) \sum_{i=1}^{N-1}
\tilde{u}_{x_i}(\bar x,t)
 \di\frac{(y_i-x_i)}{\ve}\, dy+ O(1)\chi_{\Omega_\ep}+ O(\ve^{\alpha} ).
\end{aligned}
$$

Let
$$
\begin{aligned}
B_\ep(x,t):=&\int_{\Oc}J_\ve (x-y)\Big(C_2g(\bar x,t) -
\frac{y_N-x_N}{\ve} g(\bar x,t)\Big)\\
& - \int_{\Oc} J_\ve (x-y) \sum_{i=1}^{N-1}
\tilde{u}_{x_i}(\bar x,t)
 \di\frac{(y_i-x_i)}{\ve}\, dy.
\end{aligned}
$$
Proceeding in a similar way as in the proof of Lemma \ref{F} we
get for $\ep$ small,
$$
\begin{aligned}
& \int_{\Oc}J_\ve (x-y)\Big(C_2g(\bar x,t) -
\frac{y_N-x_N}{\ve} g(\bar x,t)\Big)\\
& = g(\bar x,t)\int_{(\Oc)\cap\{|y_N-\bar x_N|
\leq \kappa
\ep^{2} \}}J_\ve
(x-y)\left( C_2 - \di\frac{(y_N-x_N)}\ep\right)\,dy\\
& \qquad + g(\bar x,t)\int_{(\Oc)\cap\{y_N-\bar x_N
> 0 \}}J_\ve
(x-y)\left( C_2 -  \di\frac{(y_N-x_N)}\ep\right)\,dy\\
& \qquad - g(\bar x,t)\int_{(\Oc)\cap\{0<y_N-\bar x_N < \kappa
\ep^{2} \}}J_\ve
(x-y)\left( C_2 - \di\frac{(y_N-x_N)}\ep\right)\,dy\\
& =C_1 g(\bar x,t)\int_{\{z_N>s\}} J(z)( C_2 - z_N)\,dz +
O(\ep^{})\chi_{\Omega_\ep}.
\end{aligned}
$$
And
$$
\begin{aligned}
&  \int_{\Oc} J_\ve (x-y) \sum_{i=1}^{N-1}
\tilde{u}_{x_i}(\bar x,t)
 \di\frac{(y_i-x_i)}{\ve}\, dy\\
&= \sum_{i=1}^{N-1} \tilde{u}_{x_i}(\bar x,t)\int_{\{|y_N-\bar x_N
|\le \kappa\ep^{2} \}}J_\ve (x-y)
 \di\frac{(y_i-x_i)}{\ve} dy\\
& \qquad + \sum_{i=1}^{N-1} \tilde{u}_{x_i}(\bar
x,t)\int_{\{y_N-\bar x_N
> \kappa\ep^{2} \}}J_\ve (x-y)
 \di\frac{(y_i-x_i)}{\ve} dy\\
& = C_1\sum_{i=1}^{N-1} \tilde{u}_{x_i}(\bar
x,t)\int_{\{z_N-s>\kappa\ep^{}  \}}J (z)z_i
dz+O(\ep^{})\chi_{\Omega_\ep}\\
 &= I_2+O(\ep^{})\chi_{\Omega_\ep}.
\end{aligned}
$$
As in Lemma \ref{F} we have $I_2=0$. Therefore,
$$
B_\ep(x,t)=C_1g(\bar x,t)\int_{\{z_N>s\}}J(z)\big(C_2 -
z_N\big)\,dz+ O(\ep)\chi_{\Omega_\ep}.
$$

Now, we observe that $B_\ve$ is bounded and supported in
$\Omega_\ep$. Hence
$$
\int_0^t
\int_\Omega |F_\ve (x,\tau)| \, dx \, d\tau \le
\frac{1}{\ve} \int_0^t \int_{\Omega_\ve} |B_\ve (x,\tau)| \, dx \, d\tau
+ Ct |\Omega_\ep| + C t |\Omega| \ve^{\alpha} \le C.
$$
This proves the first assertion of the lemma.

Now, let us write for a point $x\in\Omega_\ep$
$$
x=\bar x- \mu\eta (\bar x) \qquad\mbox{with } 0<\mu<d\ep.
$$
For $\ep$ small and $0<\mu<d\ep$, let $dS_\mu$ be the area element
of $\{x\in\Omega\,|\, \mbox{dist}\,(x,\partial\Omega)=\mu\}$.
Then, $dS_\mu=dS+O(\ep)$, where $dS$ is the area element of
$\partial\Omega$.

So that, taking now $\mu=s\ep$ we get for any continuous test
function $\theta$,
$$\begin{aligned} &\frac1\ep\int_0^T\int_{\Omega_\ep}
B_\ep(x,t)\theta
(\bar x,t)\,dx\,dt\\
&\quad =O(\ep)+{C_1}\int_0^T\int_{\partial\Omega}g(\bar x,t)
\theta (\bar x,t) \int_0^d\int_{\{z_N>s\}}J(z)\big(C_2
-  z_N\big)
\,dz\,ds\, dS\,dt \\
&\quad \quad = O(\ep)\to0\quad\mbox{as }\ve \to 0,
\end{aligned}
$$
since we have chosen $C_2$ so that
$$
\int_0^d\int_{\{z_N>s\}}J(z)\big(C_2 -  z_N\big)\,dz\,ds=0.
$$

Now, with all these estimates, we go back to $F_\ep$. We have
$$
F_\ep (x,t) = \frac{1}{\ep} B_\ep (x,t) + O (1) \chi_{\Omega_\ep}
+ O (\ep^\alpha).
$$
Thus, we obtain
$$
\int_0^T\int_{\Omega_\ep}
F_\ep(x,t)\theta (\bar x,t)\,dx\,dt\to0\qquad\mbox{as }\ep\to 0.
$$

Now, if $\sigma(r)$ is the modulus of continuity of $\theta$,
$$
\begin{aligned}
&\int_0^T\int_{\Omega_\ep} F_\ep(x,t)\theta
(x,t)\,dx\,dt=\int_0^T\int_{\Omega_\ep} F_\ep(x,t)\theta (\bar
x,t)\,dx\,dt \\
& \quad +\int_0^T\int_{\Omega_\ep}
F_\ep(x,t)\big(\theta(x,t)-\theta
(\bar x,t)\big)\,dx\,dt\\
& \le \int_0^T\int_{\Omega_\ep} F_\ep(x,t)\theta (\bar
x,t)\,dx\,dt + C\sigma(\ep)\int_0^T\int_{\Omega_\ep}
|F_\ep(x,t)|\,dx\,dt \to 0 \qquad\mbox{as }\ep\to0.
\end{aligned}
$$

Finally, the observation that $F_\ep=O(\ep^\alpha)$ in
$\Omega\setminus\Omega_\ep$ gives
$$
\int_0^T\int_{\Omega\setminus\Omega_\ep}F_\ep(x,t)\theta
(x,t)\,dx\,dt\to 0 \qquad\mbox{as }\ep\to0
$$
and this ends the proof.
\end{proof}

Now we prove that $u^\ve$ is uniformly bounded when $G=G_2$.

\begin{lema}\label{lema.que.diga} Let $G= G_2$.
There exists a constant $C$ independent of $\ve$ such that
$$
\|u^\ve \|_{L^\infty (\overline{\Omega} \times [0,T])} \le C.
$$
\end{lema}

\begin{proof}
Again we will use a comparison argument. Let us look for a
supersolution. Pick an auxiliary function $v$ as a solution to
\begin{equation}\label{calor.8}
\begin{cases}
v_t-\Delta v=h(x,t)\quad&\mbox{in}\quad \Omega\times(0,T),\\[10pt]
\di\frac{\partial v}{\partial\eta}=g_1(x,t) \quad&\mbox{on}\quad
\partial\Omega\times(0,T),\\[10pt]
v(x,0)=v_1(x)\quad&\mbox{in}\quad \Omega.
\end{cases}
\end{equation}
for some smooth functions $h(x,t)\geq 1$,  $v_1(x)\geq u_0 (x)$
and
$$
g_1 (x,t)\ge \frac2K (C_2 +1) \max_{\partial \Omega \times [0,T]}
|g(x,t)| +1
\qquad (K
\mbox{ as in }\eqref{integrales.II})
$$
such that the resulting $v$ has an extension $\tilde v$ that
belongs to $C^{2+\alpha, 1+\alpha/2}(\mathbb{R}^N\times [0,T])$
and let $M$ be an upper bound for $v$ in $\bar
\Omega\times[0,T]$.
As before $v$ is a solution to
$$
\begin{cases}
v_t-L_\ep v=H(x,t,\ep)\quad&\mbox{in}\quad \Omega\times(0,T),\\
v(x,0)= v_1(x)\quad&\mbox{in}\quad \Omega,
\end{cases}
$$
where $H$ verifies
$$
H(x,t,\ep)
\geq  \Big(\frac{ g_1(\bar x,t) }{\ve} \int_{\Oc} J_\ve (x-y) \eta(\bar x)\cdot
\frac{(y-x)}{\ve}\,dy -D_1\int_{\Oc}J_\ve (x-y)\, dy
\Big)+\frac{1}{2}.
$$
So that, by Lemma \ref{lomitos},
$$
H(x,t,\ep)
\geq \big(\frac{g_1(\bar x,t)\,K}\ep-D_1\big)\int_{\Oc}J_\ve (x-y)\, dy +\frac12
$$
for $\ep<\bar\ep$.

Let us recall that
$$
\begin{aligned}
F_\ve (x,t) =  \displaystyle
\tilde{L}_\ve (\tilde u) - \Delta \tilde u
&+\frac{C_2}{\ve}\int_{\Oc} J_\ve (x-y) g (y,t)\, dy \\
&\displaystyle -\frac{1}{\ve^2}\int_{\Oc}J_\ve
(x-y)\big(\tilde{u}(y,t)-\tilde{u}(x,t)\big)\, dy.
\end{aligned}
$$
Then, proceeding once again as in Lemma \ref{F} we have,
$$
\begin{aligned}
|F_\ep(x,t)|&\le \frac{|g (\bar x,t)| \, C_2}{\ve}\int_{\Oc} J_\ve
(x-y)
\, dy +
\frac{|g (\bar x,t)|}{\ep}\int_{\Oc} J_\ep(x-y)\big|\eta(\bar x)\cdot \frac{(y-x)}{\ve}\big| \,dy\\
& \qquad \quad + C\ep^{\alpha}+ C \int_{\Oc}J_\ve (x-y)\,dy \\
& \le \Big[\frac{(C_2+1)}\ep
\max_{\partial \Omega \times [0,T]} |g(x,t)|+C\Big]\int_{\Oc}J_\ve (x-y)\,dy +C\ep^{\alpha}\\
& \le \big(\frac {g_1(\bar x,t)\,K}{2\ep}+C\big)\int_{\Oc}J_\ve
(x-y)\,dy + C\ep^{\alpha}
\end{aligned}
$$
if $\ep<\bar\ep$, by our choice of $g_1$.

Therefore, for every $\ve$ small enough, we obtain
$$
|F_\ve (x,t)| \le H(x,t, \ve),
$$
and, by a comparison argument, we conclude that
$$
-M \le -v(x,t) \le u^\ve (x,t) \le v(x,t) \le M,
$$
for every $(x,t)\in \overline{\Omega} \times [0,T]$. This ends the
proof.
\end{proof}

Finally, we prove our last result, Theorem \ref{teo.conv.J}.

\begin{proof}[Proof of Theorem \ref{teo.conv.J}]
By Lemma \ref{lema.3.2.II} we have that
 $$ F_\ve (x,t) \rightharpoonup 0\quad\mbox{as measures in}\quad
\Omega \times [0,T]
$$
as $\ve \to 0$.

Assume first that $\psi \in C^{2+\alpha}_0 (\Omega )$ and let
$\tilde\varphi_\ep$ be the solution to
$$
\begin{aligned}
w_t-L_\ep w&= 0\\
w(x,0)&= \psi (x).
\end{aligned}
$$

Let $\tilde\varphi$ be a solution to
$$
\left\{
\begin{array}{l}
\varphi_t-\Delta \varphi = 0\\[10pt]
\displaystyle \frac{\partial \varphi}{\partial \eta} = 0 \\[10pt]
\varphi (x,0) = \psi (x).
\end{array}
\right.
$$

Then, by Theorem \ref{zero} we know that $\tilde \varphi_\ep \to
\tilde\varphi$ uniformly in $\Omega \times [0,T]$.

For a fixed $t > 0$ set $\varphi_\ve (x,s)= \tilde\varphi_\ve
(x,t-s)$. Then $\varphi_\ve$ satisfies
$$
\begin{aligned}
&\varphi_s+L_\ep \varphi =0, \qquad \mbox{ for } s<t,\\
&\varphi(x,t)=\psi (x).
\end{aligned}
$$

Analogously, set $\varphi (x,s)= \tilde\varphi (x,t-s)$. Then
$\varphi $ satisfies
$$
\left\{
\begin{array}{l}
\varphi_t+\Delta \varphi = 0\\[10pt]
\displaystyle \frac{\partial \varphi}{\partial \eta} = 0 \\[10pt]
\varphi(x,t) = \psi (x).
\end{array}
\right.
$$

Then, for $w^\ve = u^\ve -u$ we have
$$
\begin{aligned}
&\displaystyle \int_\Omega w^\ve (x,t) \,\psi (x) \, dx = \int_0^t
\int_\Omega \frac{\partial w^\ve}{\partial s} (x,s)\, \varphi_\ve (x,s) \,dx\,ds+
\int_0^t \int_\Omega \frac{\partial \varphi_\ve}{\partial s} (x,s)\,
w^\ve (x,s)\,dx\,ds
\\
&\hskip.3cm\displaystyle  = \int_0^t \int_\Omega
 L_\ve (w^\ve )(x,s) \varphi_\ve (x,s) \, dx\,
 ds
  + \int_0^t \int_\Omega F_\ve (x,s)\,
\varphi_\ve (x,s) \, dx \, ds
\\
&\hskip.8cm\displaystyle  + \int_0^t \int_\Omega
\frac{\partial \varphi_\ve}{\partial s} (x,s)\, w_\ve (x,s)\,dx\,ds\\
&\hskip.3cm =  \int_0^t
\int_\Omega L_\ve (\varphi_\ve) (x,s) w^\ve (x,s) \,  dx\,
 ds
+ \int_0^t \int_\Omega F_\ve (x,s)\,
\varphi_\ve (x,s) \, dx \, ds
\\
&\hskip.8cm\displaystyle   + \int_0^t \int_\Omega
\frac{\partial \varphi_\ve}{\partial s} (x,s)\, w^\ve (x,s)\,dx\,ds \\
 &\hskip.3cm= \int_0^t \int_\Omega F_\ve (x,s)
\varphi_\ve (x,s) \, dx \, ds.
\end{aligned}
$$
Now we observe that, by the Lemma \ref{lema.3.2.II},
$$
\begin{aligned}
& \left| \int_0^t \int_\Omega F_\ve (x,s)
\varphi_\ve (x,s) \, dx \, ds \right| \le \left| \int_0^t \int_\Omega F_\ve (x,s)
\varphi(x,s) \, dx \, ds\right| \\
& \hskip2cm +  \displaystyle \sup_{0<s<t}\|\varphi_\ve
(x,s)-\varphi(x,s)\|_{L^\infty(\Omega)}
\int_0^t \int_\Omega |F_\ve (x,s)| \, dx \, ds\to 0
\end{aligned}
$$
as $\ep\to 0$. This proves the result when $\psi\in
C^{2+\alpha}_0(\Omega)$.

Now we deal with the general case. Let $\psi \in L^1(\Omega )$.
Choose $\psi_n
\in C^{2+\alpha}_0(\Omega)$ such that $\psi_n \to \psi$ in
$L^1(\Omega)$. We have
$$
\left| \int_\Omega w^\ve (x,t) \,\psi (x) \, dx \right| \le
\left| \int_\Omega w^\ve (x,t) \, \psi_n (x) \, dx \right| + \|
\psi_n - \psi \|_{L^1 (\Omega)} \| w^\ve \|_{L^\infty (\Omega)}.
$$
By Lemma \ref{lema.que.diga}, $\{w^\ve\}$ is uniformly bounded,
and hence the result follows.
\end{proof}

\bigskip

\noindent{\bf Acknowledgements.} Supported by Universidad de Buenos Aires
under grants X052 and X066, by ANPCyT PICT No. 03-13719, Fundacion
Antorchas Project 13900-5, by CONICET (Argentina) and by FONDECYT
Project 1030798 and Coop. Int. 7050118 (Chile).

\end{document}